\documentclass[mathpazo]{cicp}
\usepackage{graphicx}

\usepackage{verbatim}

\usepackage{amssymb}
\usepackage{amsmath} 
\usepackage{latexsym}
\usepackage{subfig}
\usepackage{caption}
\usepackage{hyperref}
\usepackage{soul}
\usepackage{floatrow} 
\usepackage{url}
\usepackage{xcolor}
\usepackage{tikz} 
\usetikzlibrary{fit,shapes.geometric}

\newcommand{\added}[1]{\textcolor{black}{#1}}

\begin{document}
\title{Deep learning-based reduced-order methods for fast transient dynamics}
\author{
  Martina Cracco$^1$, Giovanni Stabile$^{1,3}$, Andrea Lario$^1$, Armin Sheidani$^{1,4}$, \\
  Martin Larcher$^2$, Folco Casadei$^2$, Georgios Valsamos$^2$, \\
  Gianluigi Rozza$^1${\corrauth}
}
\address{
  $^1$SISSA, International School for Advanced Studies, Mathematics Area, mathLab, via Bonomea 265, 34136 Trieste, Italy \\
  $^2$European Commission, Joint Research Centre, Via Enrico Fermi 2749, 21027 Ispra (VA), Italy \\
  $^3$The Biorobotics Institute, Sant'Anna School of Advanced Studies, V.le R. Piaggio 34, 56025, Pontedera, Pisa, Italy\\
  $^4$Fluids and Flows group, Department of Applied Physics and Science Education, Eindhoven University of Technology, Eindhoven, P.O. Box 513, 5600 MB, The Netherlands}
\email{{\tt grozza@sissa.it}}
\begin{abstract}
In recent years, large-scale numerical simulations played an essential role in estimating the effects of explosion events in urban environments, for the purpose of ensuring the security and safety of cities. Such simulations are computationally expensive and, often, the time taken for one single computation is large and does not permit parametric studies. The aim of this work is therefore to facilitate real-time and multi-query calculations by employing a non-intrusive Reduced Order Method (ROM).

We propose a deep learning-based (DL) ROM scheme able to deal with fast transient dynamics. In the case of blast waves, the parametrised PDEs are time-dependent and non-linear. For such problems, the Proper Orthogonal Decomposition (POD), which relies on a linear superposition of modes, cannot approximate the solutions efficiently. The piecewise POD-DL scheme developed here is a local ROM based on time-domain partitioning and a first dimensionality reduction obtained through the POD. Autoencoders are used as a second and non-linear dimensionality reduction. The latent space obtained is then reconstructed from the time and parameter space through deep forward neural networks.
The proposed scheme is applied to an example consisting of a blast wave propagating in air and impacting on the outside of a building. The efficiency of the deep learning-based ROM in approximating the time-dependent pressure field is shown.
\end{abstract}
\maketitle

\section{Introduction}
 Blast waves in cities can be caused by a terrorist attack or an accidental explosion and can provoke human casualties, different types of injuries and damage on buildings and other infrastructure. Due to the emerging threat for urban environments concerning terrorist attacks in the last decades, computational simulations have been increasingly used to identify vulnerabilities and propose protective solutions for modern cities. Recent investigations include the design of access control points (Larcher \textit{et al.}~\cite{Larcher2018}), the risk assessment in a transport infrastructure (Valsamos \textit{et al.}~\cite{Valsamos2019}) and the large-scale simulation of the explosion occurred in the port of Beirut in 2020 (Valsamos \textit{et al.}~\cite{Valsamos2021}).

 For explosions in an urban environment, one has to distinguish between the direct effects if the explosive is attached or close to a structure whose material is heavily damaged by the blast and the far-field effects of the so-called blast wave that is the focus of this work. Typically a free blast wave propagates spherically from the source of the explosion. Its pressure characteristic can be calculated in the purely spherical (or hemispherical) case by several empirical equations (Karlos \textit{et al.}.~\cite{Karlos2016}). In general, it can be said that the pressure wave is decreasing by a cubic root with the distance. But the spherical propagation might be greatly altered by reflections, shadowing and chanelling. A pressure wave hitting a rigid structure is reflected and its pressure intensity is increased by a substantial amount. Similarly, structures can shadow the area behind them since the pressure waves cannot reach these areas directly. In that case, the pressure might be much smaller and might also not be critical anymore. In the case of explosions in tube- or channel-like environments, the decrease of the pressure wave due to spherical propagation does not occur and the pressure wave amplitude might remain similar also at much bigger distances. 
 The effects of such blast waves in the far-field concern mainly the consequences on humans (Solomos \textit{et al.}~\cite{Solomos2020}) and the loading of lightweight structures, in particular windows. The use of laminated glass (Larcher \textit{et al.}~\cite{Larcher2012}) can reduce the risk of humans inside the buildings dramatically. Standards for testing such kind of structures under different blast loading are available and currently under review (Larcher \textit{et al.}~\cite{Larcher2016}). 

In industrial applications, many problems require dealing with real-time or multi-query scenarios. In the case of blast waves in urban environments caused by an explosion, being able to perform parametric studies and estimate the damage to buildings and the risk to humans in real time is crucial. For instance, in the disaster intervention units, it is necessary to understand the damage level on critical infrastructures after a blast event and adapt the intervention emergency strategy accordingly. Since these sorts of problems are to be solved in a usually very large domain (an entire city in some specific cases), they are characterised by a high-dimensional discrete space. Hence, the solution by means of a full-order method (FOM) remains prohibitively expensive in  real-time and many-query contexts. Reduced Order Methods (ROM) aim at reducing the computational time employed by a single simulation by reducing the dimension of the system.
Reduced Basis (RB) methods are examples of reduced order methods where the solutions are obtained in an offline-online computational fashion (Hesthaven \textit{et al.}~\cite{Hesthaven2015}). In the offline stage, a reduced basis space is constructed from a set of snapshots computed at different times and parameter values using a FOM solver. An approximation to the full-order solution for new time and parameter values is reconstructed in the online stage as a linear combination of the RB functions.

In this work, we present a data-driven ROM scheme able to approximate the evolution of a blast wave. An equation-free approach is preferable in cases where, such as this one, the full-order solver is not easily accessible. The ROM is based on a set of high-fidelity solutions computed using EUROPLEXUS (Casadei \textit{et al.}~\cite{casadei2001transient}), a software developed jointly by the European Commission's Joint Research Centre (JRC) and by the French Commissariat {\`a} l'{\'E}nergie Atomique et aux {\'E}nergies Alternatives (CEA), for the simulation of fast transient phenomena involving fluid and structure interaction (\cite{EUROPLEXUS2016}).

\begin{figure}[t]
\centering
\includegraphics[width = .45\textwidth,trim={0.0cm .3cm .0cm .2cm},clip]{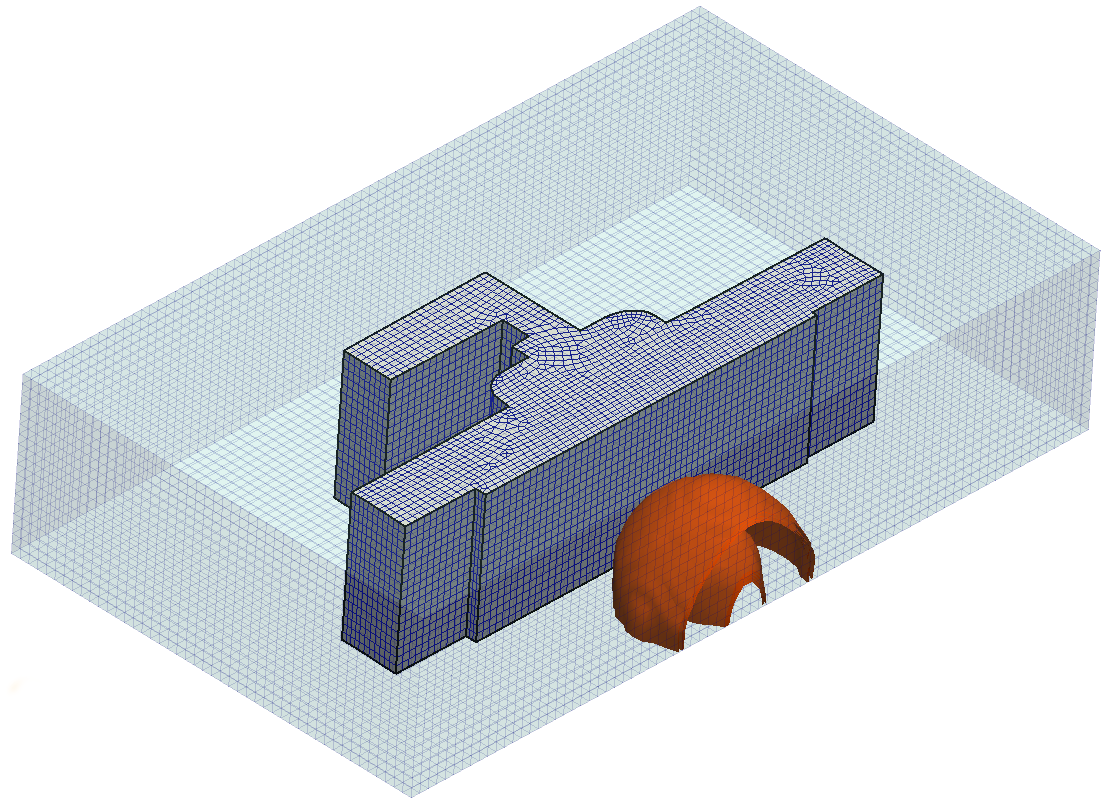}
\includegraphics[width = .45\textwidth,trim={0.0cm .3cm .0cm .2cm},clip]{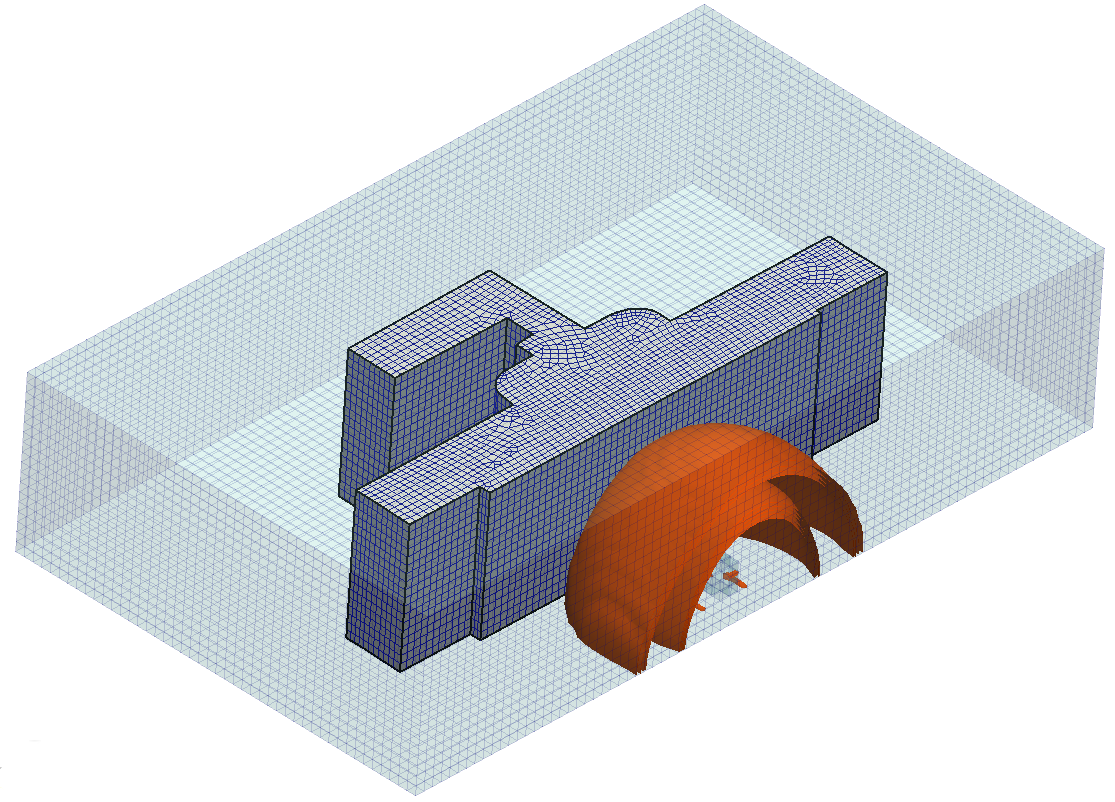}
\caption{Blast wave propagating on the outside of a building.}
\label{fig:sissa3d}
\end{figure}

In the next section, after a brief review of the literature on non-intrusive ROMs and deep learning-based ROMs, we describe the proposed reduced order scheme. In Section~\ref{sec:blast_waves}, the ROM scheme is applied to an example of a blast wave propagating in the vicinity of a building (Figure~\ref{fig:sissa3d}) and a comparison is made with the traditional POD method. Section~\ref{sec:conclusions} is dedicated to concluding remarks. 

\section{Methodology}
\label{sec:methodology}

In the case of blast waves, the parametrised PDEs are time-dependent and
non-linear and represent a transient and fast event. For such problems, the
POD which relies on a linear superposition of
modes, cannot approximate the solutions efficiently. The POD has been
successfully employed in combination with Gaussian Process Regression (GPR) or
Neural Networks (NN) to approximate the coefficients of the reduced model (e.g.\
Hesthaven and Ubbiali~\cite{Hesthaven2018}, Guo and Hesthaven~\cite{Guo2019}).
Recent works include the application of a POD-NN method to a combustion problem
(Wang \textit{et al.}~\cite{Wang2019}), a hybrid projection/data-driven ROM to
turbulent flows (Hijazi \textit{et al.}~\cite{Hijazi2020}, Georgaka \textit{et
al.} \cite{GeorgakaStabileStarRozzaBluck2020}) the combination of DMD and POD with
interpolation to shape optimisation (Tezzele \textit{et
 al.}\cite{Tezzele2019,TezzeleDemoStabileMolaRozza2020}). Recently, Lario \textit{et al.}~\cite{Lario2021}
 proposed a method which combines the spectral POD with a recurrent neural
 network to reconstruct the time evolution of the reduced coefficients.

 A non-linear alternative to linear ROMs consists of local ROMs, where an approximation is sought by partitioning the solution space into sub-regions and by assigning to each sub-region a local reduced-order basis. For example, Drohmann \textit{et al.}~\cite{Drohmann2011} and 
 Dihlmann \textit{et al.}~\cite{Dihlmann2012} used adaptive time-domain partitioning to build local ROMs.
 Local ROMs were developed further by Amsallem~\cite{Amsallem2012}, who employed an unsupervised learning algorithm in order to cluster snapshots. Amsallem \textit{et al.}~\cite{Amsallem2015} applied local ROMs in the context of hyper-reduction and Amsallem and Haasdonk~\cite{Amsallem2016} proposed a state-space partitioning criterion based on the true projection error instead of the Euclidean distance.

In recent years, several contributions have shown the potential of Deep Learning (DL) techniques also in the reduction stage of reduced order models. For example, Wiewel \textit{et al.}~\cite{Wiewel2019} proposed an approach based on Convolutional Neural Networks (CNN) and Long Short-Term Memory (LSTM) to predict the temporal evolution of a physical function. Physics-Informed Neural Networks (PINN) and Physics-Reinforced Neural Networks (PRNN) have been successfully used in the data-driven solution and discovery of partial differential equations  (e.g.\ Raissi \textit{et al.}~\cite{Raissi2019}, Chen \textit{et al.}~\cite{Chen2020}).

Another possible nonlinear and deep learning-based approach in model order reduction is based on Autoencoders. These are unsupervised artificial neural networks (ANN) able to learn a reduced representation of the given data, commonly used in image recognition and denoising (e.g.\ Gondara~\cite{Gondara2016}, Luo \textit{et al.}~\cite{Luo2018}). Autoencoders have also been successfully employed in model order reduction (e.g.\ Nikolopoulus  \textit{et al.}~\cite{Nikolopoulos2021}, Maulik \textit{et al.}~\cite{Maulik2021}, Romor \textit{et al.}~\cite{RomorStabileRozza2022}, Ahmed  \textit{et al.}~\cite{Ahmed2021}).
Fresca \textit{et al.}~\cite{Fresca2021} proposed a DL-ROM scheme, based on convolutional autoencoders, applied to PDEs representing wave-like behaviour. Applying an Autoencoder directly to the snapshots would lead to an input layer with as many neurons as degrees of freedom of the problem. This approach is feasible whenever the discrete space is not too large, and the training of ANNs in the reduction step is not excessively time-expensive. To overcome this difficulty, Phillips \textit{et al.}~\cite{Phillips2020} and Fresca \textit{et al.}~\cite{Fresca2022} employ a POD as the first reduction step, and subsequently an Autoencoder to further reduce the dimension of the problem.  \added{POD for reduced-order modelling of transient dynamics by virtue of isogeometric analysis
was developed in \cite{RichenLi2021,Zhu2017,Zhu2016}}.

To the best of our knowledge, not much work has been done on reduced order methods for blasting. For example, Xiao~\textit{et al.}~\cite{Xiao2017} used a POD-RBF scheme for solids interacting with compressible fluid flows focusing on crack propagation. 

\begin{figure}[t]
\centering
\begin{tikzpicture}
\node[shape=circle, draw=black] (uh) at (-1,0) {$\mathbf{u}_h$};
\node[shape=circle, draw=black] (u_N) at (3,0) {$\mathbf{u}_N$};
\node[shape=circle, draw=black] (u_n) at (3,-2) {$\mathbf{u}_n$};
\node[shape=rectangle, draw=black] (param_J) at (-1,-2) {$(t,\boldsymbol{\mu}), t\in T_j$};
\node[shape=rectangle, draw=black] (param) at (-5.5,-2) {$(t,\boldsymbol{\mu})$};
\path[thick,->] (param) edge node[above,swap] {\footnotesize{time-partition}} (param_J);
\path[thick,->] (uh) edge node[above,swap] {\footnotesize{POD}} (u_N);
\path[thick,->] (u_N) edge node[right,swap] {\footnotesize{$f_{\mathrm{E}}$}} (u_n);
\path[thick,->] (param_J) edge node[below,swap] {\footnotesize{DFNN}} (u_n);
\path[thick,->] (param_J) edge node[left,swap] {\footnotesize{FOM}} (uh);
\end{tikzpicture}
\caption{Piecewise POD-DL scheme: offline stage. It should be noted that $f_E$ is a function which maps the bases of the POD onto the AE laetnt space}
\label{fig:scheme}
\end{figure}
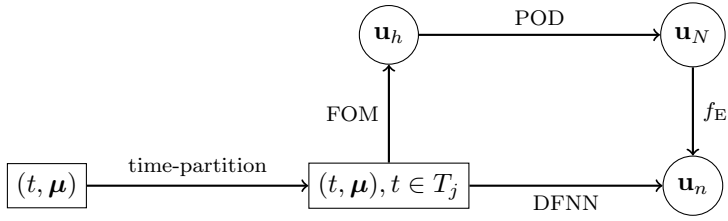

In this work, we introduce a piecewise POD-DL scheme which combines a local POD approach and deep learning methods able to deal with fast and transient phenomena involving fluid-structure interactions such as blast waves in the vicinity of buildings. Local ROMs based on time-domain partitioning are particularly suited for problems characterised by different physical regimes and fast development, such as blast waves. We apply a first reduction layer which consists of a piecewise-POD in time, followed by an Autoencoder (AE) which works as a second non-linear reduction step. Finally, we train a Deep Forward Neural Network (DFNN) to learn the latent space dynamics. The offline procedure is summarised in Figure~\ref{fig:scheme}.
The novelty of this work lies in the use of deep learning techniques combined with a local approach and the application to a large-scale system representing fast and transient dynamics. It should be noted that the code repository associated with this paper is available at: \url{https://github.com/Mardgi/sissa-jrc}. \added{It is important to note that, while this method has significant potential for broad applications in fluid mechanics, the focus of this study, as a collaboration between SISSA and JRC, was on investigating the code for fast transient dynamics due to the associated challenges with generalizing.}
\subsection{Proper Orthogonal Decomposition}
\label{sec:pod}
In this section, we revise the Proper Orthogonal Decomposition method (Benner \textit{et al.}\cite{book2021}). Let $u_h$ be a high-fidelity solution to a parametrised partial differential equation that can be written as follows 
\[u_h (t,\mathbf{ x}, \boldsymbol{\mu}) = \sum_{i=1}^{N_h} u_h^{(i)} (t,\boldsymbol{\mu})\*\phi_i(\mathbf{ x}), \]
where $\{\phi_i\}$ represents a base of the discrete approximation space $V_h$ (eg. finite volume method), $N_h = \dim(V_h) $ is the number of degrees of freedom, $t \in \mathcal{T}$ (time interval), $\mathbf{ x} \in \Omega \subset \mathbb{R}^3$ (space domain) and $\boldsymbol{\mu} \in  \mathcal{P} \subset \mathbb{R}^d$ (parameter domain). The vector $\mathbf{u}_h = [u_h^{(1)},  \dots, u_h^{(N_h)}] \in \mathbb{R}^{N_h}$ collects the coefficients of the approximation.

The time domain $\mathcal{T}$ is partitioned into $J$ sub-intervals $T_j$ such that $\mathcal{T} = \cup_{j=1}^{J} T_j$. Let $\mathcal{P}_h = \{\boldsymbol{\mu}_j\}_{j=1}^{N_\mu}$ be a sampling of the parameter space $\mathcal{P}$. We can obtain a collection of high-fidelity solutions by applying a full-order solver for different parameter values in $\mathcal{P}_h$ and divide them according to the time-domain partitioning. For each subset of snapshots, we apply the same procedure. For the sub-interval $T_j$, we build the snapshot matrix $S$ defined as follows
\[S = [S_1, \dots, S_{N_{\boldsymbol{\mu}}}] \in \mathbb{R}^{N_h \times N_s},\]
where $N_s = N_{\boldsymbol{\mu}} \* N_t$ is the number of snapshots and $S_k$ is the time-trajectory matrix that collects solutions for the parameter value $\boldsymbol{\mu}_k$ and  $t_i \in T_j$, i.e.
\[S_k = [\mathbf{u}_h (t_1, \boldsymbol{\mu}_k), \dots, \mathbf{u}_h (t_{N_t}, \boldsymbol{\mu}_k)] \in \mathbb{R}^{N_h \times N_t}.\]

We seek a reduced basis approximation to the high-fidelity solution $\mathbf{u}_h(t,\boldsymbol{\mu})$ for $t \in T_j$ in the form
\begin{equation}
\mathbf{u}_{h,\mathrm{POD}}(t,\boldsymbol{\mu}) = \sum_{i=1}^N {u}^{(i)}_{N}(t,\boldsymbol{\mu})\*\mathbf{v}_i = V\*\mathbf{u}_{N},\label{eq:pod}
\end{equation}
such that $N \ll N_s$. The vector $\mathbf{u}_{N} \in \mathbb{R}^N$ collects the reduced coefficients, while $\{\mathbf{v}_i\}$ form the basis of a reduced space $V_\mathrm{rb}$.

We can use the well-known POD method in order to extract an orthogonal basis. We apply a singular value decomposition to the snapshot matrix $S$, which takes the form
\[S = U\*\Sigma\*W^T,\] 
where $U, W$ are unitary matrices and $\Sigma  = \text{diag}(\sigma_1, \dots, \sigma_r, 0, \dots, 0)$ is a diagonal matrix with decreasing singular values on the diagonal and rank $r \leq \min (N_h, N_s)$. 

The basis consisting of the first $N<r$ left singular vectors of $S$ is the best approximation in the least-squares sense to the snapshots collected in $S$. The orthogonal matrix  $V$, composed of the first $N$ column vectors of $U$,  minimises the projection error of the snapshots over all orthogonal matrices in $\mathbb{R}^{N_h \times N}$, i.e.
\[\|S- V\*V^T\*S\|_F =\min_{\substack{Z \in \mathbb{R}^{N_h\times N}\\ Z^T\* Z = \mathbb{I}_I}} \|S- Z\*Z^T\*S\|_F = \sum_{i=N+1}^r \sigma_i^2,\]
with $\|\cdot\|_F$ being the Frobenius norm and $\mathbb{I}_N$ the identity matrix.
Therefore, the amount of information carried by the first $N$ modes depends on the square of the singular values of the system. We can set a criterion to choose the dimension of the basis $N = \dim (V_{rb})$ such that
\[\frac{\sum^r_{i=N+1} \sigma_i^2}{\sum^r_{i=1} \sigma_i^2} \leq \epsilon,\]
for a chosen tolerance $\epsilon$. The time partitioning is determined so that the projection error in each time interval is of the same order. 

\subsection{Autoencoders}
\added{An Autoencoder is a neural network architecture consisting of two interconnected parts, designed so that the input and output have the same dimensions. The first part, called the \textit{encoder}, maps an input vector $\mathbf{z} \in \mathbb{R}^N$ to a lower-dimensional representation $f_\mathrm{E}(\mathbf{z}) \in \mathbb{R}^n$, where $n < N$. The second part, known as the \textit{decoder}, takes this compressed representation $f_\mathrm{E}(\mathbf{z}) \in \mathbb{R}^n$ and reconstructs it to match the input dimensions, i.e., $f_\mathrm{D}(f_\mathrm{E}(\mathbf{z})) \in \mathbb{R}^N$. Autoencoders are considered unsupervised learning methods because they learn patterns without the need for labeled data.
}

In this work, we exploit autoencoders in order to reduce further after the application of the POD as described in Section~\ref{sec:pod}. Specifically, we train an autoencoder to compress the vector of reduced coefficients $\mathbf{u}_N=V^T\*S \in \mathbb{R}^N$ in the POD approximation defined by equation~\eqref{eq:pod}, to a latent space representation $\mathbf{u}_n \in \mathbb{R}^n$ with $n<N$. Once the AE is trained, we have an encoder map $f_\mathrm{E}$ and a decoder map $f_\mathrm{D}$ that can be used to find the reconstructed coefficients $\mathbf{u}_{N,\mathrm{AE}}$ as follows
\[\mathbf{u}_{N,\mathrm{AE}} = f_\mathrm{D}(f_\mathrm{E}(\mathbf{u}_N)) \in \mathbb{R}^N.\]
An approximation to the full-order solution $\mathbf{u}_h$ for $t \in T_j$ can be obtained as a linear combination of the reduced basis functions, i.e.
\begin{equation}
     \mathbf{u}_{h,\mathrm{AE}} = V\*\mathbf{u}_{N,\mathrm{AE}} \in \mathbb{R}^{N_h}.\label{eq:ae}
 \end{equation}

\subsection{Regression}

\begin{figure}[t]
\centering
\begin{tikzpicture}
\node[shape=ellipse, draw=black] (uh) at (6,-2) {$\mathbf{u}_{h,\mathrm{POD-DL}}$};
\node[shape=circle, draw=black] (u_N) at (6,0) {$\mathbf{u}_N$};
\node[shape=circle, draw=black] (u_n) at (3.5,0) {$\mathbf{u}_{n}$};
\node[shape=rectangle, draw=black] (param_j) at (0,0) {$(t,\boldsymbol{\mu}), t\in T_j$};
\node[shape=rectangle, draw=black] (param) at (0,-2) {$(t,\boldsymbol{\mu})$};
\path[thick,->] (param) edge node[left,swap] {\footnotesize{time-switch}} (param_j);
\path[thick,->] (u_N) edge node[right,swap] {\footnotesize{POD}} (uh);
\path[thick,->] (u_n) edge node[above,swap] {\footnotesize{$f_{\mathrm{D}}$}} (u_N);
\path[thick,->] (param_j) edge node[above,swap] {\footnotesize{DFNN}} (u_n);
\end{tikzpicture}
\caption{Piecewise POD-DL scheme: online stage.}
\label{fig:scheme_online}
\end{figure}
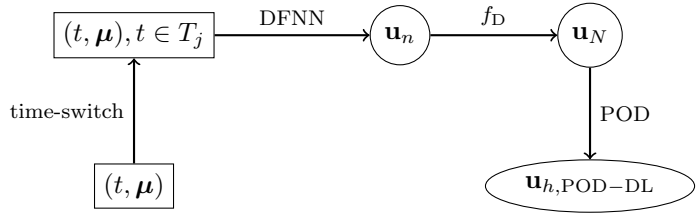

Finally, a regression model is used to approximate the map $\pi: \mathcal{T} \times \mathcal{P} \rightarrow \mathbb{R}^n$ from the time and parameter space to the latent space, defined as
\[(t, \boldsymbol{\mu}) \to \pi(t,\boldsymbol{\mu}) = \mathbf{u}_{n}(t, \boldsymbol{\mu}), \]
where $\mathbf{u}_{n}= f_\mathrm{E}(\mathbf{u}_N)$ is the encoded representation of $\mathbf{u}_N$. In this work, we use a DFNN, although other choices are possible (e.g.\ Guo and Hesthaven~\cite{Guo2019}). Once an approximated regression map $\pi_{\mathrm{NN}}$ is obtained, for any untrained time and parameter vector $(\bar{t}, \bar{\boldsymbol{\mu}})$ the approximated representation of the
solution in the latent space is obtained by evaluating the regression model, i.e.
\[\mathbf{u}_{n}(\bar{t}, \bar{\boldsymbol{\mu}})\approx {\mathbf{u}}_{n,\mathrm{NN}}(\bar{t}, \bar{\boldsymbol{\mu}}) = {\pi}_{\mathrm{NN}}(\bar{t}, \bar{\boldsymbol{\mu}}). \]
In the online stage, as illustrated in Figure~\ref{fig:scheme_online}, we recover an approximation of $\mathbf{u}_h$ as a linear combination of the POD modes, where the coefficients are obtained by applying the decoder function, $f_\mathrm{D}$ to the approximated latent space, as follows
\begin{equation}
     \mathbf{u}_{h,\mathrm{POD-DL}}(\bar{t}, \bar{\boldsymbol{\mu}}) = V\*f_\mathrm{D}(\mathbf{u}_{n,\mathrm{NN}}(\bar{t}, \bar{\boldsymbol{\mu}})) \in \mathbb{R}^{N_h}.\label{eq:uhpoddl}
\end{equation}

\subsection{Error calculations}
We can evaluate the performance of the entire scheme by calculating the error between the high-fidelity solution $\mathbf{u}_{h}$ and its approximation ${\mathbf{u}}_{h,\mathrm{POD-DL}}$, which satisfies the following inequality
\[\|\mathbf{u}_{h} - {\mathbf{u}}_{h,\mathrm{POD-DL}}  \| \le \|\mathbf{u}_{h} - {\mathbf{u}}_{h,\mathrm{POD}}  \|  + \|\mathbf{u}_{h,\mathrm{POD}} - {\mathbf{u}}_{h,\mathrm{AE}}  \|  +\|\mathbf{u}_{h,\mathrm{AE}} - {\mathbf{u}}_{h,\mathrm{POD-DL}}  \|,  \]
where $\|\cdot\|$ is the chosen norm in $\mathbb{R}^{N_h}$. In this work, we consider the $L^1,L^2$ and $L^\infty$ vector norms. Thus, we define the relative projection error
\begin{equation}\epsilon_{\mathrm{POD}}(t, \boldsymbol{\mu}) = \frac{\|\mathbf{u}_h - \mathbf{u}_{h,\mathrm{POD}}  \|}{\|\mathbf{u}_h\|} = \frac{\|\mathbf{u}_h - V\*V^T\*\mathbf{u}_h  \|}{\|\mathbf{u}_h\|},\label{eq:error_pod}\end{equation}
with $\mathbf{u}_{h,\mathrm{POD}}$ being the projection onto the reduced space defined by equation~\eqref{eq:pod}. To evaluate the performance of the Autoencoder alone, we define the relative error
\begin{equation}\epsilon_{\mathrm{AE}}(t, \boldsymbol{\mu}) = \frac{\|\mathbf{u}_{h,\mathrm{POD}} - {\mathbf{u}}_{h,\mathrm{AE}}  \|}{\|\mathbf{u}_{h,\mathrm{POD}}\|}, \label{eq:error_ae}\end{equation}
where ${\mathbf{u}}_{h,\mathrm{AE}}$ is defined by equation~\eqref{eq:ae}. The relative regression error, which accounts only for the error committed by the DFNN, is defined as
\begin{equation}\epsilon_{\mathrm{NN}}(t, \boldsymbol{\mu}) = \frac{\|\mathbf{u}_{h,\mathrm{AE}} - {\mathbf{u}}_{h,\mathrm{POD-DL}}  \|}{\|\mathbf{u}_{h,\mathrm{AE}}\|}.\label{eq:error_nn}\end{equation}
Finally, we evaluate the entire POD-DL scheme by calculating the relative error
\begin{equation}\epsilon_{\mathrm{POD-DL}}(t, \boldsymbol{\mu}) =\frac{\|\mathbf{u}_{h} - {\mathbf{u}}_{h,\mathrm{POD-DL}}  \|}{\|\mathbf{u}_{h} \|}. \label{eq:error_poddl}\end{equation}

\section{Blast waves in urban environments}
\label{sec:blast_waves}
The example analysed in this work consists of an explosion occurring on the outside of a building (Figure~\ref{fig:sissa3d}). We obtain a set of snapshots by computing the time-dependent pressure field and by varying the $x$-position of the explosion. 

The high-fidelity solutions are obtained with EUROPLEXUS adopting a compressed bubble model (phenomenological model) (Larcher \textit{et al.}~\cite{Larcher2010}) for the representation of the blast load in the fluid. The explosion is simulated by imposing an initial condition consisting of a bubble of high pressure, which expands and propagates throughout the fluid domain. The overpressure in the bubble is automatically calculated given the volume and the mass of the charge.

The air in the fluid domain is modelled as a perfect gas and discretised using a uniform \added{mesh of hexahedral}, cells \added{having} a uniform length in each direction. The \added{general Navier-Stokes equations, expressing the conservation of mass, momentum and energy, are reduced to their} compressible and inviscid \added{form} (Euler) equations \added{thanks to the dominance of pressure forces with respect to viscous forces and} are solved by marching in time with an explicit second-order scheme \added{under the Courant-Friedrichs-Lewy (CFL) stability condition. No turbulence model is included, as usual in this class of applications}. The building is represented by a non-deformable structure embedded in the fluid mesh and a fast search algorithm is employed to determine which fluid cells lie inside the structure's influence domain \added{, where Fluid-Structure Interaction (FSI) is modelled by an immersed body algorithm} (Casadei \textit{et al.}~\cite{Casadei2011}). The fluid is described by an Eulerian formulation, while the structure is described by a Lagrangian formulation. The \added{spatial} discretisation technique for the fluid domain is a second-order cell-centred finite volume \added{(FV)} method, where all the variables are \added{collocated} at the centre of the volumes. \added{In the Eulerian formulation used for the fluid domain, the mesh (referential domain) is fixed while the fluid moves relative to these fixed grid points. The numerical fluxes across FV cell interfaces are calculated using the approximate Harten-Lax-van-Leer-Contact (HLLC) Riemann solver\cite{toro2009riemann}. In order to achieve second-order accuracy in space, the code uses the Green-Gauss reconstruction of the conservative variables at the interfaces between fluid volumes, from the values at the volume centroids, by using the spatial gradients at the centroids and by applying an appropriate flux limiter.} The fluid domain consists of 66240 finite volumes with a mesh size of 2.5m. \added{To prevent spurious blast wave reflections, infinite} boundary conditions are applied to the surfaces of the fluid domain which interface with the open space. 
EUROPLEXUS is able to perform a risk analysis, which is based on the calculation of two important quantities, impulse and peak overpressure.
Impulse and peak overpressure at the final time are essential quantities to the estimation of the consequences on structures and humans (Ferradas \textit{et al.}~\cite{Ferradas2008}). 
Let $p_0$ denote the atmospheric pressure ($p_0 = 10^5$ Pa). Impulse and peak overpressure are defined, respectively, as
 \begin{align}
    I(t) &= \int_0^{t} \big(p(\tau)-p_0\big)_+\ d\tau,\label{eq:impulse}\\ 
    p_\mathrm{max}(t) &= \max_{\tau \in (0,t)} \big(p(\tau)-p_0\big)_+, \label{eq:overp}
 \end{align}
where $p$ is the pressure and $(\cdot)_+$ represents the positive part of a function. 

\subsection{Problem description}
For the purpose of applying the reduced order scheme introduced in Section~\ref{sec:methodology}, we generate a set of $N_\mu = 30$ pressure trajectories, i.e. $\{\mathbf{u}_h(t_j, \mu_k)\}_{j,k}$, where  $j \in \{1,\dots,N_t\}$ and $k \in \{1,\dots,N_\mu\}$. Since the high-fidelity solver works with International System of Units (SI), the pressure is expressed in Pascal while distances and time are measured in meters and seconds, respectively. From this point on, for the sake of brevity, we omit units in our results. The parameter $\mu=x_\mathrm{exp}$ represents the $x$-position of the explosion and the parameter space is $\mathcal{P} = [129, 189]$. \added{It is important to highlight that the range of parameters was chosen to cover a meaningful spatial domain where the explosion's effect on the structure is significant while avoiding regions where boundary effects might dominate the solution. This range ensures that the explosion remains within the fluid domain and generates complex pressure trajectories relevant to the study. Additionally, this range allowed us to sample a diverse set of solutions for training the reduced-order model while maintaining computational feasibility. Also, We sampled the parameter space using an equispaced approach.} \added{In this article we included examples with a one dimensional parameter because the actual parametrization already generates a rather complex solution manifold. However, from a theoretical standpoint, the methodology can also deal with parameter spaces larger than one.} The space domain where the structure is submersed is $\mathcal{D} = [60, 258]\times[130,245]\times[0,45] \subset \mathbb{R}^3$. Other fixed parameters consist of the mass of the explosive, $m_{\mathrm{exp}} = 100$, the $y$- and $z$-coordinates of the explosion, $y_{\mathrm{exp}} = 140$ and $z_{\mathrm{exp}} = 5.5$. The high-fidelity solution $\mathbf{u}_h \in \mathbb{R}^{N_h}$ represents the pressure field with $N_h = 66240$ degrees of freedom. Since a cell-centred finite volume scheme is employed, each component of the vector $\mathbf{u}_h$ represents the pressure at the centre of a cell of the fluid grid. For this problem, the number of time steps employed by the high-fidelity solver is $N_t =520$. As time domain, we restrict our attention to $\mathcal{T} = [0.01, 0.5]$, since the air bubble for $t<0.01$ is still confined to the vicinity of the initial explosion and far from the structure. The solution at the first instants is more difficult to approximate given the high overpressure concentrated in only a few elements of fluid. 
The time-interval $\mathcal{T}$ is partitioned empirically, ensuring that the projection error for each sub-region is of the same order of magnitude for a particular choice of the number of modes $N$.

The encoder part of the autoencoder is composed of 7 layers with a constant decrement of neurons between each layer. For example, for $N=200$ and $n=20$, the decoder's layers are $200\rightarrow170\rightarrow140\rightarrow110\rightarrow80\rightarrow50\rightarrow20$. The decoder network consists of the same layers in the reversed order. The Leaky Relu (Xu et al.\cite{xu2020reluplex} and Dubey et al. \cite{dubey2019comparative}) with $\alpha=0.01$ is chosen as the activation function for each layer except for the last decoding layer, where a Sigmoid function is used to ensure an output between 0 and 1. The weights are initialised using the He Normal initialiser and the optimisation algorithm is Adam Optimiser with learning rate $\lambda = 10^{-4}$. As loss function, the Mean Squared Error (MSE) is used and the Mean Absolute Error (MAE) is also monitored. An Early Stopping mechanism is used, where the training stops if the loss function on a validation set does not improve over  200 epochs. The maximum number of epochs is set to 10000 and the chosen batch size is 32.

Regarding the DFNN used as a regressor from the time and parameter space to the latent space, we use 8 layers of 50 neurons each, with Sigmoid activation functions and He Normal kernel initialiser. As for the AE, the MSE is chosen as the loss function and the MAE is also monitored. The optimisation algorithm is Adam (Zhang\cite{zhang2018improved}, Bock and Wei\ss \cite{bock2019proof}) Optimiser with learning rate $\lambda = 10^{-4}$, the maximum number of epochs is set to 5000, and the batch size to 32. The early stopping mechanism is used with a patience of 200 epochs.

The snapshots are divided into train, validation and test set with, respectively, $N_{\mathrm{tr}} = 0.75\*N_s$, $N_{\mathrm{val}} = 0.15\*N_s$ and $N_{\mathrm{te}} = 0.10\*N_s$ number of samples. The training set is used to calculate the POD modes and train the neural networks, while the validation set is used for the Early Stopping of the neural network training. The test set is used to evaluate the scheme on untrained data. \added{The computational cost of training is comparable to that of a relatively small neural network, typically completed within hours. Testing involves inference using the same network, making its computational demand negligible, on the order of milliseconds.}

\subsection{Results}

\begin{figure}[t]
\centering
\includegraphics[width=0.95\textwidth,trim={0cm .3cm 0cm .3cm},clip]{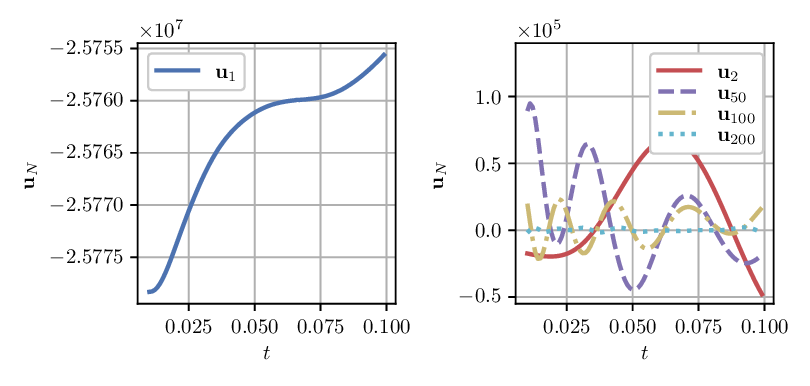}
\caption{Reduced basis coefficients for $t \in T_1$, for reduced dimension $N=200$ and latent dimension $n=20$. Coefficient relative to the 1\textsuperscript{st} mode (left), and coefficients relative to the 2\textsuperscript{nd}, 50\textsuperscript{th}, 100\textsuperscript{th} and 200\textsuperscript{th} modes (right).}
\label{fig:reduced_coeffs}
\end{figure}
\begin{figure}
\centering
\includegraphics[trim={0cm .3cm 0cm .4cm},clip]{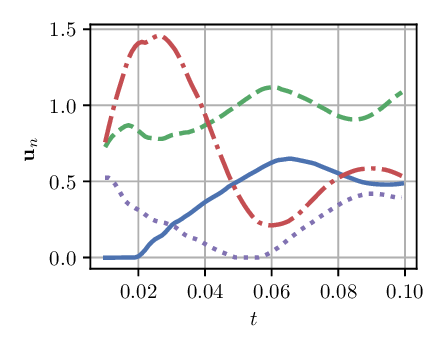}
\caption{Time evolution of the first four latent representation vectors for reduced dimension $N=200$, latent dimension $n=20$ and for $t \in T_1$.}
\label{fig:latent_space}
\end{figure}

\begin{figure}[t]
\centering
\sidesubfloat[]{\includegraphics[width = .92\textwidth,trim={0.0cm .3cm .0cm .2cm},clip]{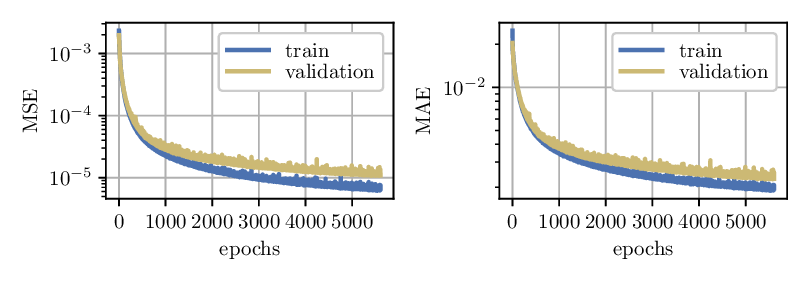}}

\sidesubfloat[]{\includegraphics[width = .92\textwidth,trim={0.0cm .3cm .0cm .2cm},clip]{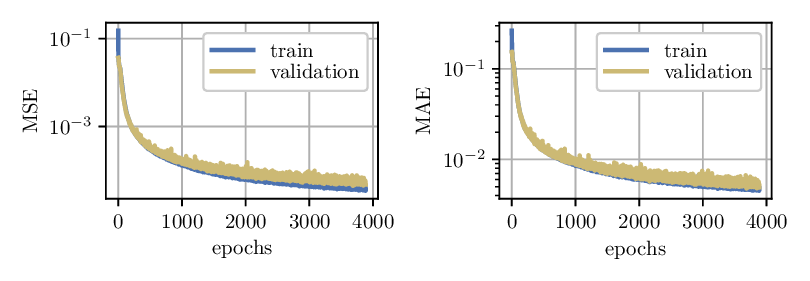}}
\caption{Training history for reduced space dimension $N=200$, latent dimension $n=20$ and $t\in T_1$. Training and validation MSE (loss function) on the left, training and validation MAE (monitored metric) on the right. (a) Autoencoder; (b) Deep Forward Neural Network.}
\label{fig:train_history}
\end{figure}

In Figure~\ref{fig:reduced_coeffs}, we show the reduced order coefficients $\mathbf{u}_N$ for $t \in T_1$ and for number of modes  $N=200$ and latent dimension $n=20$. 
We notice that the first reduced basis coefficient $\mathbf{u}_1$ is two orders of magnitude higher than the second coefficient $\mathbf{u}_2$ for the time interval $T_1$. Moreover, the following coefficients decrease in magnitude and show an increasingly oscillatory behaviour. For this reason, we need to apply a suitable normalisation before training the autoencoder to ensure that all coefficients be accurately approximated. Since the first coefficient $\mathbf{u}_1$, relative to the dominant mode, is negative we decide to normalise it between 0 and 1 separately from the rest of the coefficients, as follows
\[\overline{\mathbf{u}}_1 = \frac{\mathbf{u}_1 - \min{(\mathbf{u}_1)}}{\max{(\mathbf{u}_1)} - \min{(\mathbf{u}_1)}},\]
where the minimum and maximum are to be intended as calculated over all times in the training set. The remaining coefficients are  normalised between 0 and 1, as follows
\[\overline{\mathbf{u}}_N = \frac{\mathbf{u}(2:) - \min{(\mathbf{u}(2:))}}{\max{(\mathbf{u}(2:))} - \min{(\mathbf{u}(2:))}},\]
where the minimum and maximum are calculated over all times and for all $N>2$ and the notation ($u(2:)$) indicates the second mode onwards.

The time evolution of the first four latent vectors is shown in Figure~\ref{fig:latent_space} for $t \in T_1$, for reduced dimension $N=200$ and latent space dimension $n=20$ \added{at $z=5$ plane}. \added{It should be noted that the z=5 plane corresponds to the region of the domain that is closest to the explosion source in the 
z-direction. At this depth, the pressure wavefront interacts strongly with the structure, making this plane the most critical for evaluating the solution's accuracy. By focusing on this plane, we assess the reduced-order model's performance in the region where the physical phenomena of interest are most pronounced.} 
The latent space vectors $\mathbf{u}_n \in \mathbb{R}^n$ are normalised between 0 and 1 before training the DFNN, as follows
\[\overline{\mathbf{u}}_n = \frac{\mathbf{u}_n - \min{(\mathbf{u}_n)}}{\max{(\mathbf{u}_n)} - \min{(\mathbf{u}_n)}},\]
where the minimum and maximum are over all times and for all $n$.

In Figure~\ref{fig:train_history}, the training history for autoencoder and deep forward neural network are shown in the case where $N=200, n=20$ and $t \in T_1$. In Figure~\ref{fig:train_history}(a), we show the loss function, MSE, and the MAE in the training and validation sets for the AE network. Notice that the early stopping mechanism determines the stop of the optimisation algorithm at about 5700 epochs. In  Figure~\ref{fig:train_history}(b), the MSE and MAE calculated in the training and validation sets for the DFNN are represented. We can observe that the optimisation stops at about 3900 epochs.

\begin{figure}[p]
\centering
\sidesubfloat[]{\includegraphics[width=.9\textwidth,trim={.0cm .0cm .0cm .0cm},clip]{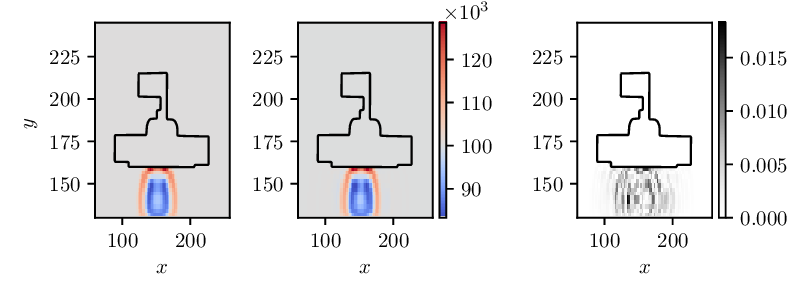}}

\sidesubfloat[]{\includegraphics[width=.9\textwidth,trim={.0cm .0cm .0cm .0cm},clip]{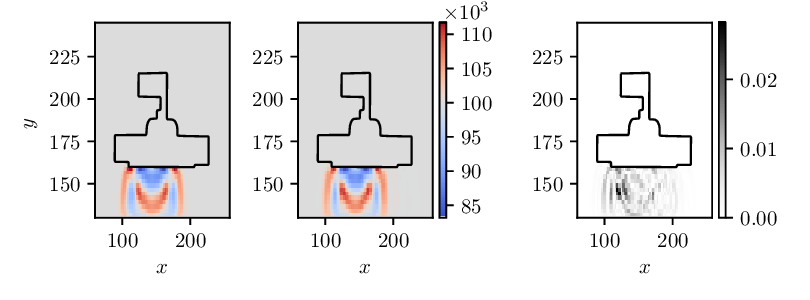}}

\sidesubfloat[]{\includegraphics[width=.9\textwidth,trim={.0cm .0cm .0cm .0cm},clip]{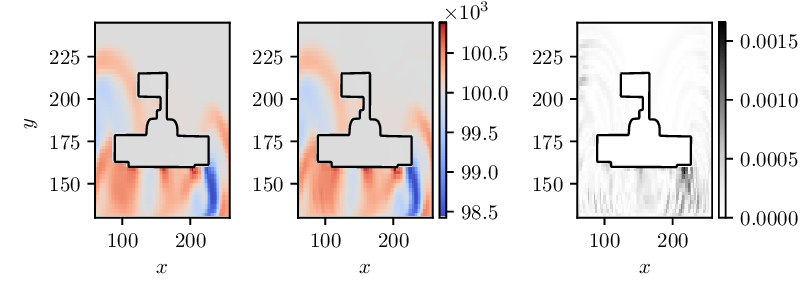}}

\sidesubfloat[]{\includegraphics[width=.9\textwidth,trim={.0cm .0cm .0cm .0cm},clip]{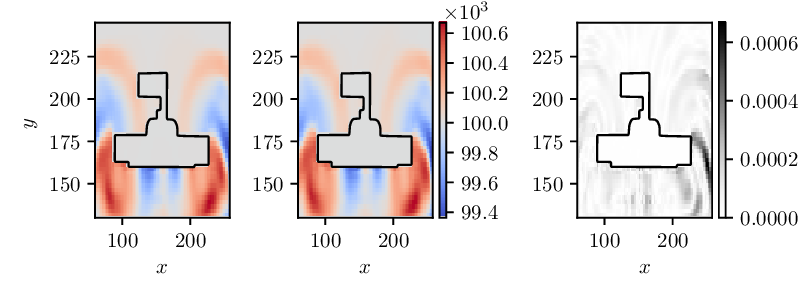}}

\caption{Comparison between high-fidelity pressure field, $\mathbf{u}_h$ (left), and POD-DL solution,  $\mathbf{u}_{h,\mathrm{POD-DL}}$ (centre), using $N=200$, $n=20$. Relative error on the right. At (a) $(\bar{t}, \bar{x}_\mathrm{exp}) = (0.07, 170.38)$; (b) $(\bar{t}, \bar{x}_\mathrm{exp}) = (0.1, 131.07)$; (c) $(\bar{t}, \bar{x}_\mathrm{exp}) = (0.2, 160.03)$; (d) $(\bar{t}, \bar{x}_\mathrm{exp}) = (0.4, 180.72)$.}
\label{fig:pressure}
\end{figure}

Figure~\ref{fig:pressure} shows a comparison of the high-fidelity pressure field $\mathbf{u}_{h}$ and the POD-DL solution $\mathbf{u}_{\mathrm{POD-DL}}$ at $z=5$ for different times and explosion positions $(\bar{t},\bar{x}_\mathrm{exp})$ in the test set. The relative error is calculated as 
\begin{equation}\epsilon(x,y) = \frac{|\mathbf{u}_{h}- \mathbf{u}_{\mathrm{POD-DL}} |}{\mathrm{mean}{(\mathbf{u}}_{h})},\label{eq:error_level} \end{equation} 
where the mean is over the domain at $z=5$. We notice that the POD-DL scheme is able to reconstruct the solution accurately. The errors are higher at smaller times, e.g.\ for $(\bar{t},\bar{x}_\mathrm{exp}) = (0.07,170.38)$ the relative errors reach approximately $1.5\%$ while for $(\bar{t},\bar{x}_\mathrm{exp}) = (0.4, 180.72)$ they are up to $0.06\%$. At $t\approx 0.01$, the pressure wave is still confined to a small portion of the domain, where the overpressure is extremely high with respect to the rest of the domain. In this case, the solution is more difficult to approximate even when retaining a large number of modes in the POD expansion and this leads to larger errors. 

\begin{figure}[t]
\centering
\includegraphics[width=0.95\textwidth,trim={.0cm .1cm .0cm .1cm},clip]{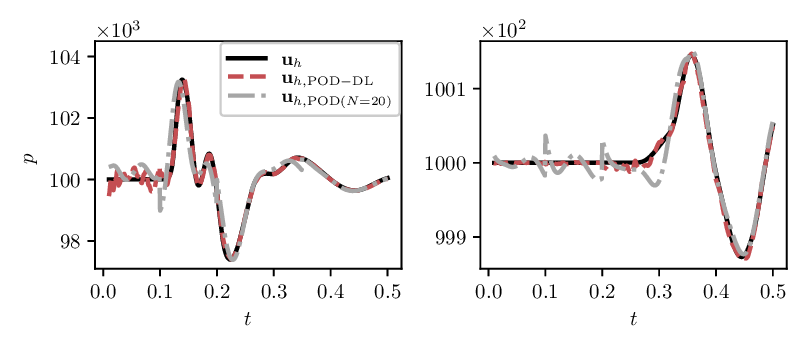}
\caption{High-fidelity pressure trajectory $\mathbf{u}_{h}$, piecewise POD solution $\mathbf{u}_{h,\mathrm{POD}}$ with $N=20$ and POD-DL trajectory $\mathbf{u}_{h,\mathrm{POD-DL}}$ obtained with $N=200$ and $n=20$ for $x_\mathrm{exp} =158.20$, at two fixed positions in the space domain. At  $(x,y,z) = (100,150,5.5)$ (left),  $L^2$-error in time for POD: 17.30, for POD-DL: 5.07. At $(x,y,z) = (200,200,5.5)$ (right),  $L^2$-error in time for POD: 0.67, for POD-DL: 0.20.}
\label{fig:pressure_dt}
\end{figure}

In Figure~\ref{fig:pressure_dt}, we present a comparison between the piecewise POD-DL scheme and the piecewise POD in reconstructing the pressure trajectory at two fixed positions in the space domain $\mathcal{D}$ for an explosion at $x_\mathrm{exp} = 158.20$. The POD-DL trajectories $\mathbf{u}_{h,\mathrm{POD-DL}}(t)$ are reconstructed by applying the reduced scheme as described by equation \eqref{eq:uhpoddl} for the same time steps used by the FOM to obtain $\mathbf{u}_{h}(t)$. Note that $\mathbf{u}_{h,\mathrm{POD}}$ represents the projected solution, i.e.\ $\mathbf{u}_{h,\mathrm{POD}}= V\*V^{T}\*\mathbf{u}_{h}$, without any regression or interpolation of the reduced coefficients. Therefore, it represents the best that can be obtained with the POD method, without accounting for the errors introduced by reconstructing the coefficients from the parameter space.
We can see that the scheme hereby proposed with $N=200$ and $n=20$ performs considerably better compared to a piecewise POD with $N=20$ modes per each time interval. For a point located near the left corner of the building $(x,y,z) = (150,100,5.5)$, the overpressure reaches about $3000$ Pa and the $L^2$-error in time amounts to 17.30 for the piecewise POD scheme and 5.07 for the piecewise POD-DL scheme. For $(x,y,z) = (200,200,5.5)$ located behind the building, the overpressure amounts to $200$ Pa and the $L^2$-error in time is equal to 0.67 for the piecewise POD scheme and to 0.20 for the piecewise POD-DL scheme.

In this particular case, the time employed by EUROPLEXUS to compute one trajectory is approximately $445 s$, while the ROM method described in Section~\ref{sec:methodology} requires about $0.2 s$ to compute an approximated solution in the online phase. It is important to notice that, the latter refers to the time required to obtain an approximated solution for given values of the time and $x-$position of the explosion. In order to obtain the time required to compute the whole trajectory, it is necessary to multiply it by the appropriate number of \added{time stations}. 

\subsubsection{Impulse and Overpressure}
\begin{figure}[p]
\RawFloats
\centering
\sidesubfloat[]{\includegraphics[width=0.85\textwidth,trim={.0cm .1cm .0cm .1cm},clip]{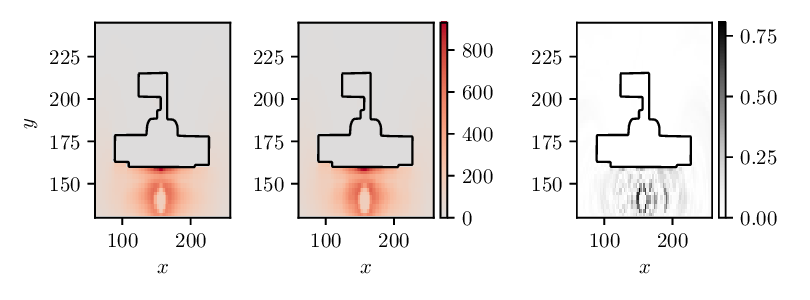}}

\sidesubfloat[]{\includegraphics[width=0.85\textwidth,trim={.0cm .1cm .0cm .1cm},clip]{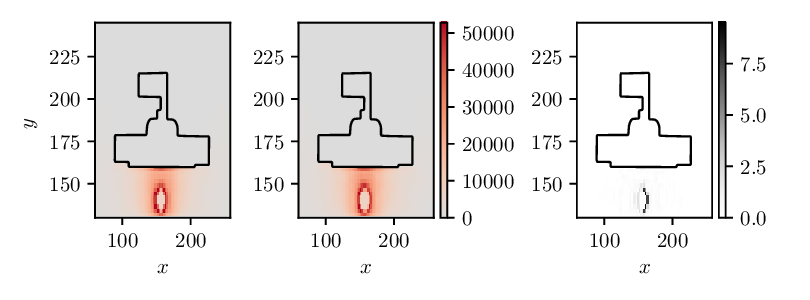}}

\caption{Comparison of the high-fidelity solution (left) and the reduced solution (centre) using $N=200$ and $n=20$ at $z=5$, for $\bar{x}_\mathrm{exp}=158.20$, and relative error (right). (a) Impulse at $t=0.5$; (b) Peak Overpressure at $t=0.5$.}
\label{fig:imp_overp}
\bigskip

\sidesubfloat[]{\includegraphics[width=0.85\textwidth,trim={.0cm .1cm .0cm .1cm},clip]{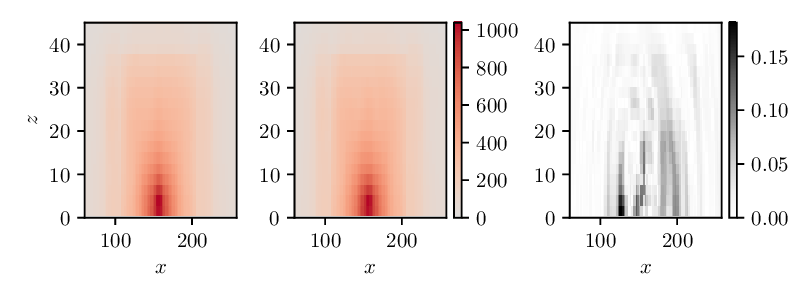}}

\sidesubfloat[]{\includegraphics[width=0.85\textwidth,trim={.0cm .1cm .0cm .1cm},clip]{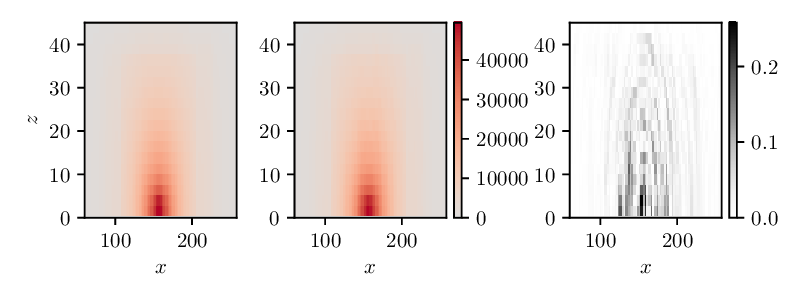}}
\caption{Comparison of the high-fidelity solution (left) and the reduced solution (centre) using $N=200$ and $n=20$ on the front facade ($y=160$) of the building, for $\bar{x}_\mathrm{exp}=158.20$, and relative error (right). (a) Impulse at $t=0.5$; (b) Peak Overpressure at $t=0.5$.}
\label{fig:imp_overp_face}
\end{figure}

In Figure~\ref{fig:imp_overp}, we show the impulse and overpressure for an explosion located at $x_\mathrm{exp}=158.20$, at the final time $t=0.5$ and  a height of $z=5$ from the ground. Impulse and overpressure are calculated from the pressure trajectory using equations~\eqref{eq:impulse}-\eqref{eq:overp}. We compare the high-fidelity solutions $\mathbf{I}_h, \mathbf{p}_{\mathrm{max},h}$ with the reduced solution $\mathbf{I}_{\mathrm{POD-DL}}, \mathbf{p}_{\mathrm{max},\mathrm{POD-DL}}$ computed with $N=200$ and $n=20$ and we calculate the error as done for the pressure field, using equation~\eqref{eq:error_level}. We can see, from Figure~\ref{fig:imp_overp}(a), that for the impulse the relative errors amount to 0.75 in some areas of the domain. In Figure~\ref{fig:imp_overp}(b), the relative errors in reconstructing the overpressure reach 7.5 in the vicinity of the original explosion.

Figure~\ref{fig:imp_overp_face} represents the impulse and overpressure at $t=0.5$ and at a fixed $y= 160$,  which corresponds to the front facade of the building. The error between high-fidelity and reduced solutions ($N=200, n=20$) is calculated using equation~\eqref{eq:error_level}, where the mean is calculated over the domain at $y= 160$. The errors are at most $15\%$ for the impulse and $20\%$ for the overpressure. 

\begin{figure}[!t]
\centering
\sidesubfloat[]{\includegraphics[width=.91\textwidth,trim={.0cm .0cm .0cm .0cm},clip]{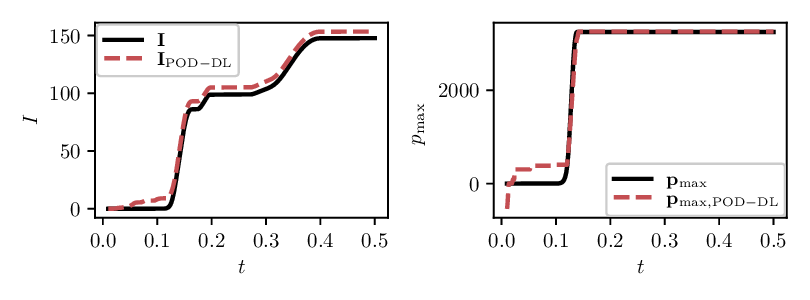}}

\sidesubfloat[]{\includegraphics[width=.91\textwidth,trim={.0cm .0cm .0cm .0cm},clip]{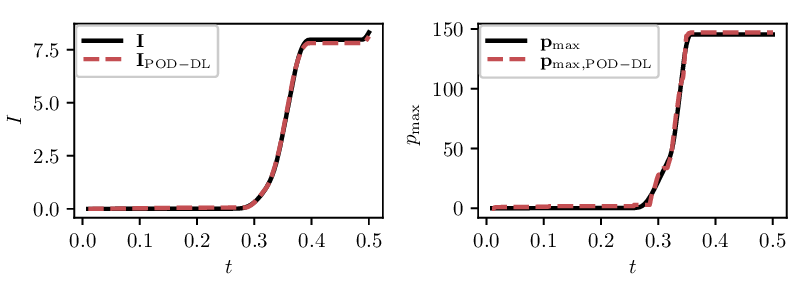}}
\caption{Comparison of the FOM and the POD-DL trajectories using $N=200$ and $n=20$, for $x_\mathrm{exp} =158.20$. (a) Impulse (left) and Peak Overpressure (right) at $(x,y,z) = (100,150,5.5)$; (b) Impulse (left) and Peak Overpressure (right) at $(x,y,z) = (200,200,5.5)$.}
\label{fig:imp_overp_dt}
\end{figure}

In Figure~\ref{fig:imp_overp_dt}, we compare the trajectories of impulse and overpressure for an explosion position $x_\mathrm{exp}=158.20$ at two different positions in the domain, $(x,y,z) = (100,150,5.5)$ located near the left corner of the building and at $(x,y,z) = (200,200,5.5)$, located behind the building.

\subsubsection{Error analysis}
\begin{figure}[t]
\RawFloats
\centering
\includegraphics[width=0.9\textwidth]{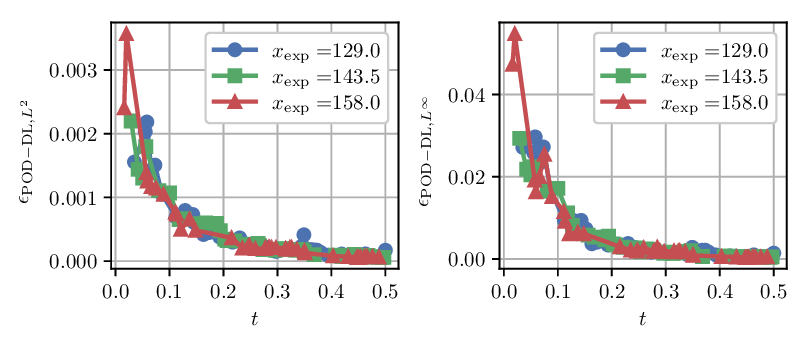}
\caption{Errors against time in $L^2$ (left) and $L^\infty$ norm (right) between the high-fidelity solution $\textbf{u}_h$ and the reduced solution $\textbf{u}_{\mathrm{POD-DL}}$ for $N=200$, $n=20$ and different values of $x_\mathrm{exp}$. }\label{fig:err_vstime}

\bigskip

\includegraphics[width=0.9\textwidth]{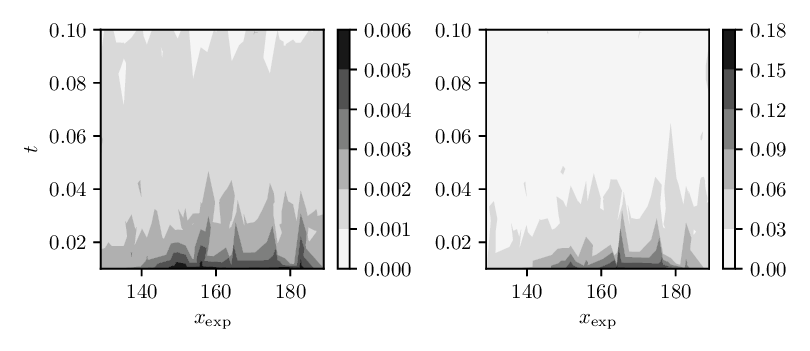}
\caption{Errors in a $x_\mathrm{exp}t$ plane between the high-fidelity solution $\textbf{u}_h$ and the reduced solution $\textbf{u}_{\mathrm{POD-DL}}$ for $N=200$ and $n=20$, using the $L^2$ norm (left) and the $L^\infty$ norm (right). }\label{fig:err_vstime_contour}
\end{figure}

In Figure~\ref{fig:err_vstime}, we present the error $\epsilon_{\mathrm{POD-DL}}$, defined by equation \eqref{eq:error_poddl}, as it varies with time and calculated in the test set for different values of $x_\mathrm{exp}$. We notice that, as already seen in Figure~\ref{fig:pressure}, the errors in $L^2$ and $L^\infty$ norm are considerably larger near $t\approx 0.01$. Moreover, the errors for $x_\mathrm{exp} =159$ are larger than the errors for $x_\mathrm{exp}=129$. This behaviour can be observed also in Figure~\ref{fig:err_vstime_contour}, where the $L^2$ and $L^\infty$ norm errors are shown in a $x_\mathrm{exp}t$ plane. The errors in the middle of the parameter space $\mathcal{P} = [129,189]$ are larger than at the boundaries of the parameter space. 

\begin{figure}[p]
\RawFloats 
\centering
\includegraphics[width=0.7\textwidth]{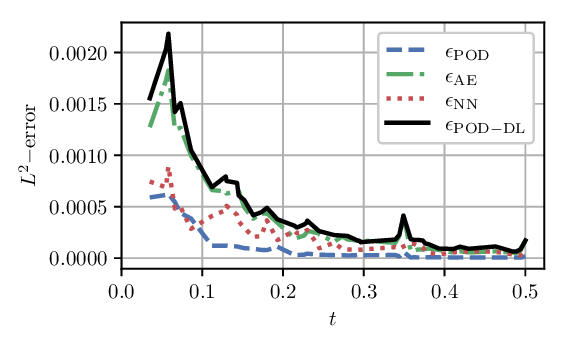}
\caption{Different components of error in $L^2$ norm against time for $N=200, n=20$ and $\bar{x}_\mathrm{exp} = 129$.}
\label{fig:err_components_vstime}

\bigskip

\centering
\includegraphics[width=1\textwidth]{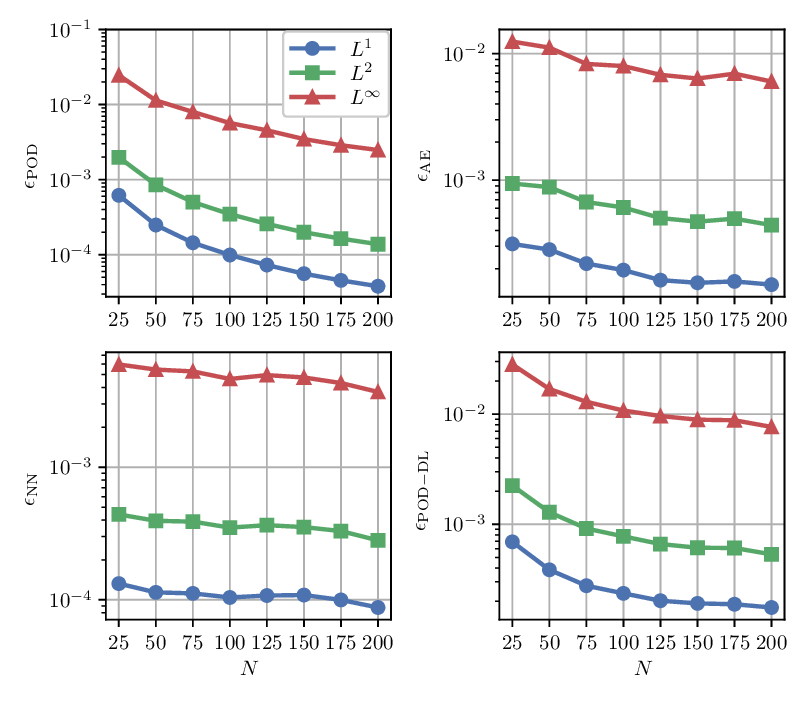}
\caption{Errors against the reduced space dimension $N$ for a fixed value of latent dimension $n=20$.}
\label{fig:err_varyN}
\end{figure}

In Figure~\ref{fig:err_components_vstime}, we split the various sources of error to look at how they contribute to the global error $\epsilon_{\mathrm{POD-DL}}$, given by equation~\eqref{eq:error_poddl}. We plot separately the POD error $\epsilon_{\mathrm{POD}}$, defined by equation~\eqref{eq:error_pod}, the error given by the autoencoder $\epsilon_{\mathrm{AE}}$, defined by equation~\eqref{eq:error_ae}, and the error given by the neural network employed as regressor $\epsilon_{\mathrm{NN}}$, defined by equation~\eqref{eq:error_nn}. We notice that all the components of the errors are larger at $t\approx0.1$ and that the autoencoder is the greatest source of error.

We perform an error analysis to investigate the performance of the piecewise POD-DL scheme, varying the number of modes $N$ retained in the POD expansion and  the latent dimension $n$ \added{as presented in Fig.\ref{fig:err_varyN}}.
First, we vary the reduced space dimension $N$ from 25 to 200, while keeping fixed the latent dimension to $n=20$. We plot the different components of error in $L^1, L^2$ and $L^\infty$ norm. We note that all the errors decay as $N$ increases. In particular, the $L^2$ POD error, $\epsilon_{\mathrm{POD}}$, is approximately $0.2\%$ for $N=25$ and decreases to about $0.01\%$ for $N=200$. The AE error in $L^2$ norm, $\epsilon_{\mathrm{AE}}$, improves slightly from $0.1\%$ for $N=25$ to about $0.05\%$ for $N=200$. As $N$ increases, there is an increasing amount of information that the autoencoder compresses to a latent space of dimension $n=20$. 
The NN error in $L^2$ norm, $\epsilon_{\mathrm{NN}}$, sees a slight improvement over $N$ since the latent space reconstructed from the time and parameter space has a fixed dimension $n=20$. 
Finally, the global error in $L^2$ norm, $\epsilon_{\mathrm{POD-DL}}$, decreases until about $0.05\%$ for $N=125$. For larger reduced dimensions, there is no significant improvement.

\begin{figure}[p]
\RawFloats
\centering
\includegraphics[width=1\textwidth]{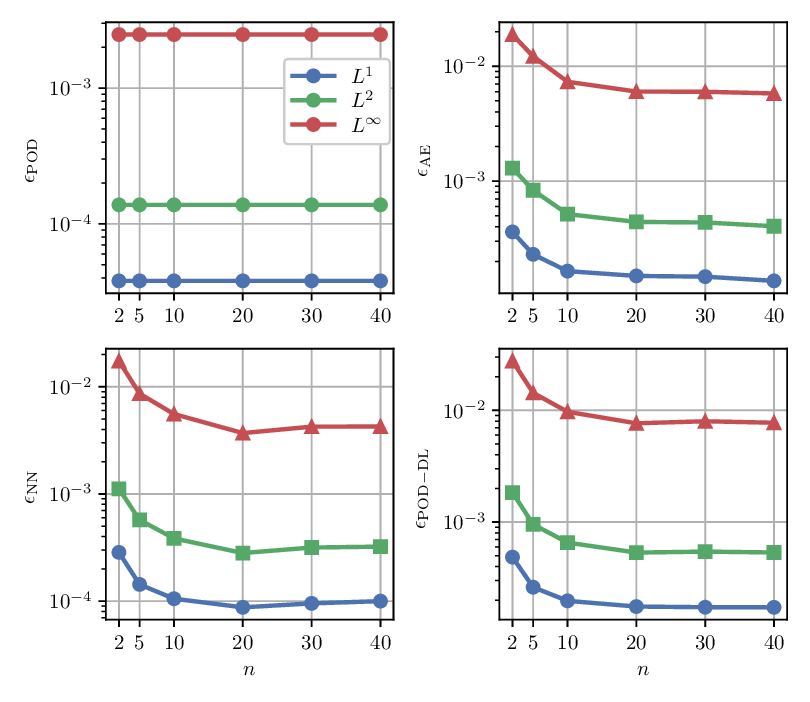}
\caption{Errors against the latent dimension $n$ for a fixed value of reduced space dimension $N=200$.}
\label{fig:err_varyn}

\bigskip

\centering
\includegraphics[width=1\textwidth]{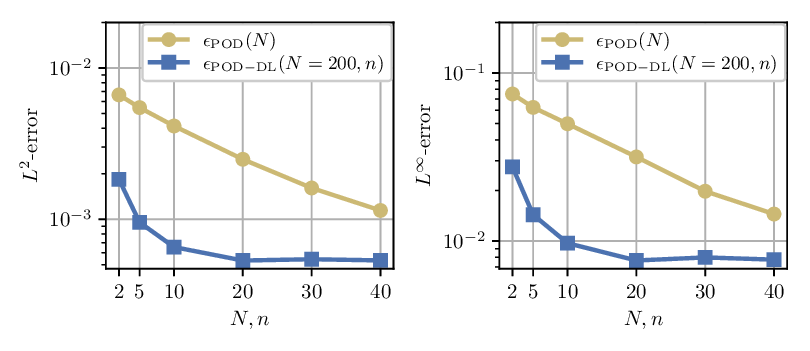}
\caption{Comparison between POD errors against $N$ and POD-DL errors against $n$ for a fixed $N=200$, in $L^2$ norm (left) and $L^\infty$ norm (right).}
\label{fig:err_vspod}
\end{figure}

Figure~\ref{fig:err_varyn} shows how errors vary with latent dimension $n$, while keeping constant the number of modes $N=200$. In this case, the errors introduced by the POD projection, $\epsilon_{\mathrm{POD}}$, remain constant. The AE errors, $\epsilon_{\mathrm{AE}}$, decay until about $n=20$ and flatten for larger values of $n$. The errors introduced by the regression, $\epsilon_{\mathrm{NN}}$, decay until about $n=20$ and then slowly increase with $n$. This is because a larger latent space is more difficult to reconstruct from the time and parameter space. The global error, $\epsilon_{\mathrm{POD-DL}}$, decreases until about $n=20$. For larger values of $n$, there is no improvement.

Figure~\ref{fig:err_vspod} shows the projection error, $\epsilon_\mathrm{POD}$, with increasing number of modes $N$ and the POD-DL  error, $\epsilon_\mathrm{POD-DL}$, with fixed $N=200$ and increasing latent dimension $n$. We can see how the POD-DL scheme leads to smaller errors in both $L^2$ and $L^\infty$ norm. In particular, $\epsilon_\mathrm{POD-DL}$ with $N=200,n=20$ is one order of magnitude smaller than $\epsilon_\mathrm{POD}$ with $N=20$.

\section{Conclusions}
\label{sec:conclusions}
In this work, we proposed a local and non-linear Reduced Order Method able to efficiently approximate the fast evolution of blast waves. This kind of problems requires the solution of non-linear and time-dependent parametric PDEs and is characterised by a large scale and fast and transient development. The validity of the method is shown in an example consisting of a blast wave propagating near a building. A set of snapshots are obtained by varying the $x$-component of the explosion position.

The scheme hereby proposed is based on a local approach, where 
the snapshots are divided following a time-domain partitioning and a local-reduced order basis is built for each sub-region. The reduced order spaces are generated using the POD and by retaining a large number $N$ of modes. The POD method alone is not able to represent the solutions efficiently since the singular values decay very slowly for this type of problems \added{by almost an order of magnitude with respect to the Autoencoder modes, which is also evident in the number of extracted mode from POD and latent space dimension}. Therefore, the remarkable compression capabilities of autoencoders are exploited to build a latent space of dimension $n < N$. A regression map between the time and parameter space and the latent space is then approximated using a deep forward neural network. 

We showed that, for new values of the time and parameter, the piecewise POD-DL method leads to very accurate approximations of the pressure field. Moreover, it performs considerably better in reconstructing the pressure trajectories compared to the POD method. Then, we tested the method in approximating the impulse and peak overpressure fields at the final time. 

Finally, we performed an error analysis. We presented the $L^2$ and $L^\infty$ norm errors as they vary with time and showed that the greatest errors are for $t\approx 0$ as expected. The solutions at smaller times consist of a high-pressure bubble confined to a small region. This causes difficulties for the POD to approximate the solutions even when retaining a high number of modes. We analysed errors by varying the number of modes $N$ and the latent dimension $n$. We proved that errors decrease when increasing the value of $N$ and keeping the latent space dimension fixed. Also, the errors decay when increasing $n$ and fixing the reduced space dimension. However, at about $n \approx 20$ the errors flatten due to the regression with DFNN. 

It is important to notice that the computational cost to recreate one particular \added{time station} with the ROM method (equal to approximately $0.2\text{s}$) would not change significantly even if a large-scale model was employed. On the contrary, the computational cost for a simulation with \added{the explicit code} EUROPLEXUS \added{increase linearly with the number of cells (for a fixed cell size) whereas with an implicit code the cost growth would be even more rapid}. That makes the proposed approach suitable for a fast (even real-time) solution once the model has been trained\added{,} even for a large-scale problem.

In terms of future perspectives, we aim to test the developed methodology also for different physics and verify the applicability of the compression strategy for the development of projection-based ROMs \cite{StabileRozza2018, RomorStabileRozza2022}. Moreover, possible improvements are in the direction of studying unsupervised learning strategies to identify clusters of solution into both the solution space, the parameter space and the time domain. 

\section*{Acknowledgements}
This work was partially funded by European Union Funding for Research and Innovation --- Horizon 2020 Program --- in the framework of European Research
Council Executive Agency: H2020 ERC CoG 2015 AROMA-CFD project 681447 ``Advanced Reduced Order Methods with Applications in Computational Fluid Dynamics'', P.I. Professor Gianluigi Rozza. 
We  acknowledge the INDAM-GNCS 2020 project ``Advanced Numerical Techniques for Industrial Applications".
We also acknowledge the support by MIUR (Italian Ministry for
University and Research) FARE-X-AROMA-CFD project and PRIN ``Numerical Analysis for Full and Reduced Order Methods for Partial Differential Equations" (NA-FROM-PDEs) project and by the FSE 2014/2020 ``Mobilità degli assegnisti di ricerca nei centri di ricerca JRC" PS 72/17. 
GS acknowledges the financial support under the National Recovery and Resilience Plan (NRRP), Mission 4, Component 2, Investment 1.1, Call for tender No. 1409 published on 14.9.2022 by the Italian Ministry of University and Research (MUR), funded by the European Union – NextGenerationEU– Project Title ROMEU – CUP P2022FEZS3 - Grant Assignment Decree No. 1379 adopted on 01/09/2023 by the Italian Ministry of Ministry of University and Research (MUR).
This work was in collaboration with the European Commission's Joint Research Centre, Unit E.4 ``Safety and Security of Buildings" in the framework of the ``Numerical Simulations of Human Brain Vulnerability to Blast Loading" project.

\cleardoublepage
\bibliographystyle{acm}
\bibliography{main}

\end{document}